\documentclass[11pt,reqno]{amsart}
\usepackage{amssymb}

\renewcommand{\le}{\leqslant}
\renewcommand{\ge}{\geqslant}
\newcommand{\RR}{\mathbb{R}}
\newcommand{\ZZ}{\mathbb{Z}}
\newcommand{\NN}{\mathbb{N}}

\newcommand{\diam}{\mathrm{diam}}

\newtheorem{theorem}{Theorem}

\newtheorem*{theoremK}{Theorem K}
\newtheorem*{theoremG}{Theorem G}
\newtheorem*{theoremKM}{Theorem KM}
\newtheorem*{theoremS}{Theorem S}
\newtheorem*{theoremVV}{Theorem VV}
\newtheorem{theoremBV}{Theorem BV}

\newcommand{\Snsim}{\mathcal{S}_n(\psi_1,\psi_2,\dots,\psi_n)}
\newcommand{\Snmult}{\mathcal{S}_n^*(\psi)}

\begin{document}

\title[Simultaneous Diophantine Approximation on Planar Curves]{Convergence Results for Simultaneous \& Multiplicative Diophantine Approximation on Planar Curves. }

\author[D.~Badziahin]{Dzmitry Badziahin}
\address{Dzmitry Badziahin, Department of Mathematics, University of
  York, Heslington, York, YO10 5DD, United Kingdom}
\email{db528@york.ac.uk}

\author[J.~Levesley]{Jason Levesley}
\address{Jason Levesley, Department of Mathematics, University of
  York, Heslington, York, YO10 5DD, United Kingdom}
\email{jl107@york.ac.uk}

\subjclass[2000]{Primary 11J83; Secondary 11J13, 11K60}

\keywords{Diophantine approximation, Khintchine type theorems,
planar curves, Hausdorff measure.}

\begin{abstract}
  Let $\mathcal{C}$ be a non-degenerate planar curve. We show that
  the curve is of Khintchine-type for convergence in the case of
  simultaneous approximation in $\RR^2$ with two independent
  approximation functions; that is if a certain sum converges then
  the set of all points $(x,y)$ on the curve which satisfy simultaneously the
  inequalities $\| q x \| < \psi_1(q)$ and $\| qy \| < \psi_2(q)$
  infinitely often has induced measure $0$. This completes the
  metric theory for the Lebesgue case. Further, for
  multiplicative approximation $\| qx \| \| q y \| < \psi(q)$ we
  establish a Hausdorff measure convergence result for the same
  class of curves, the first such result for a general class of
  manifolds in this particular setup.
\end{abstract}

\maketitle

\section{Introduction}
\label{intro}

Let $\psi:\RR^+\to \RR^+$ be a real, decreasing function. Throughout
we will refer to $\psi$ as an approximating function. Let $x\in\RR$
and $\|x\|$ be the distance of $x$ from $\ZZ$. That is, $\| x \| =
\inf\{|x -z|:z\in\ZZ\}$. Further, if $S$ is a (Lebesgue) measurable
set in $\RR^n$ then we will denote the Lebesgue measure, or more
simply the measure, of $S$ by $|S|_{\RR^n}$.

Consider now the following system of $n$ Diophantine inequalities
\begin{equation}
    \label{eq:syssimDineqs}
    \| q x_i  \| < \psi_i(q)
\end{equation}
where $x_i\in\RR$, $p_i\in\ZZ$, $q\in\NN$ and
$\psi_1,\psi_2,\dots,\psi_n$ are approximation functions. Then a
point $x\in\RR^n$ is simultaneously
$(\psi_1,\psi_2,\dots,\psi_n)$-approximable if there are infinitely
many $q$ satisfying~\eqref{eq:syssimDineqs}. The set of all such
points $x\in\RR^n$, will be denoted by
$\mathcal{S}_n(\psi_1,\psi_2,\dots,\psi_n)$.

Simultaneous approximation has another variant in the guise of
multiplicative approximation; let $x\in\RR^n$, then $x$ is
multiplicatively $\psi$-approximable if the inequality
\begin{equation}
    \label{eq:syssimmultDineqs}
    \prod_{i=1}^n{}\| q x_i \| < \psi(q)
\end{equation}
holds for infinitely many $q\in\NN$. By analogy with the previous
notion of Simultaneous approximation, we shall denote by $\Snmult$
the set of all multiplicatively $\psi$-approximable points $x$ in
$\RR^n$.

The following results, the first of which is due to Khintchine and
is a generalisation of Khintchine's own result of 1924 which deals
with the case when $\psi_1 =  \psi_2 = \dots = \psi_n $. The second
is due to Gallagher. Together they give an almost complete answer to
the question of the `size', in terms of $n$-dimensional Lebesgue
measure, of $\Snsim$ and $\Snmult$.

\begin{theoremK}\label{ksim}
    Let $\psi_1,\psi_2,\dots,\psi_n$ be approximation functions as
    defined above. Then
    \[
            |\Snsim|_{\RR^n} = \begin{cases}
                                    0 & \text{\quad{}if\quad} \sum
                                    \psi_1(h)\dots\psi_n(h) < \infty
                                    \\
                                    \text{FULL} & \text{\quad{}if\quad} \sum
                                    \psi_1(h)\dots\psi_n(h) = \infty
                               \end{cases}
    \]
\end{theoremK}

\begin{theoremG}\label{gmult}
    Let $\psi$ be an approximating function. Then
    \[
            |\Snmult|_{\RR^n} = \begin{cases}
                                    0 & \text{\quad{}if\quad} \sum
                                    \psi(h)^n\log^{n-1}h < \infty
                                    \\
                                    \text{FULL} & \text{\quad{}if\quad} \sum
                                    \psi(h)^n\log^{n-1}h = \infty
                               \end{cases}
    \]
\end{theoremG}

It is to be understood that the term `full' means that the
complement of the set in question is of measure $0$.

In Theorem(s) K and G the approximation problem is an independent
variables problem; no functional  relationship exists between any of
the coordinates of $x$. Once a functional relationship is assumed to
exist , that is the points $x$ are constrained to lie on a
sub-manifold $\mathcal{M}\subset\RR^n$, then the corresponding
approximation problems become a great deal more difficult. As we
shall see below, until very recently hardly any general results at
all were known for such problems. Furthermore, of those results that
have been established, most hold only for planar curves with
sufficient curvature conditions. There is next to nothing known for
manifolds of dimension in $\RR^n$ where $n\geq{}3$. This is in stark
contrast to the dual form of approximation; see for
example~\cite{BBKM},~\cite{BersVelSlice} or ~\cite{BKM}, where the
state of knowledge is much more complete. Before we discuss some of
the more significant results which are known to hold in the
so-called ``dependent variables''\footnote{the terminology is due to
Sprind\v{z}uk~\cite{Sprindzuk}} case, it is necessary to define the
concept of non-degeneracy.

Let $C^{(m)}(U)$ be the space of all $m$--continuously
differentiable functions $f$ where $f:U\to\RR$ with $U$ being an
open set in $\RR^n$. A map $g:U\to\RR$ is said to be non-degenerate
at $u\in{}U$ if there exists some $l\in\NN$ such that
$g\in{}C^{(l)}(B(u,\delta))$ for some sufficiently small $\delta>0$
with $B(u,\delta)\subset{}U$ and the partial derivatives of $g$
evaluated at $u$ span $\RR^n$. The map $g$ is non-degenerate if it
is non-degenerate at almost all points $u\in{}U$. Let $\mathcal{M}$
be a sub-manifold of $\RR^n$. Then $\mathcal{M}$ is said to be
non-degenerate if $\mathcal{M}=g(U)$ where $g$ is non-degenerate.
The geometric interpretation of non-degeneracy is that the manifold
is curved enough that it deviates from any hyperplane.
Non-degeneracy is not a particularly restrictive condition and a
large class of manifolds satisfy this condition.

Note that if the topological dimension, $\dim\mathcal{M}$, of the
manifold is strictly less than $n$ then $|\mathcal{M}|_{\RR^n} = 0$.
As we wish to make measure theoretic statements about points that
lie on the manifold we work with the induced measure,
$|\cdot|_{\mathcal{M}}$. All ``zero-full'' statements are made with
respect to this restricted measure.

Schmidt~\cite{S} proved one of the first major results that
attempted to address the problem of generalising Theorems K
and G 
to manifolds.
\begin{theoremS}
    Let $\psi:\RR^+\to\RR:x\mapsto{x^{-\tau}}$ where $\tau>0$ and
    $2\tau>1$. Then for any $C^{(3)}$ non-degenerate planar curve,
    $\mathcal{C}$,
    \[
        |\mathcal{C}\cap\mathcal{S}_2(\psi,\psi)|_\mathcal{C} = 0.
    \]
\end{theoremS}
Note that this result is not quite a Khintchine-type result for
convergence as the measure $0$ statement does not depend on the
convergence of an associated sum. This question was settled very
recently for an arbitrary approximation function $\psi$.
In~\cite{BDV}, Beresnevich, et. al.\ established the divergence part
of such a theorem by showing that for any $C^{(3)}$ non-degenerate
curve
\begin{equation}\label{BDVV1}
    |\mathcal{C}\cap\mathcal{S}_2(\psi,\psi)|_\mathcal{C} =
    \text{FULL\quad{}if\quad}
    \sum_{h=1}^\infty \psi^2(h) = \infty.
\end{equation}
To prove this result the authors adapted their notion of ``local
ubiquity'' as developed in~\cite{BDV0401118}. Interestingly the
convergence case initially alluded them. This is yet another
manifestation of the difficulties one encounters when trying to
establish dependent variable analogues of many of the classical
results of Diophantine approximation. In the classical setting the
convergence case is usually straightforward and all the substance is
in the divergence case. However, for dependent variable problems it
turns out that both halves of a Khintchine-type result are highly
non-trivial with the convergence case sometimes turning out to be
the more difficult of the two.

In the above case the main obstacle to establishing the convergence
statement was the need for precise information about the number of
rational points near the curve. This is a notoriously difficult
problem and the best known result, due to Huxley ~\cite[Theorem
4.2.4]{Hux}, was good enough to ensure a reasonable distribution of
rational points near the curve for the application of local
ubiquity, but gave estimates which were too large for certain sums
when considering the convergence case. Using a result of Vaughan's,
which appeared as an appendix in~\cite{BDV} and was a significant
sharpening of Huxley's result for the class of rational quadrics,
curves that are the image of the unit circle, the parabola $\{(x,y):
y=x^2\}$ or the hyperbola $\{ (x,y): x^2 - y^2 =1\}$ under a
rational affine transformation of the plane, the authors were able
to establish the convergence counterpart of this theorem for such
manifolds. Subsequently Beresnevitch \& Velani~\cite{BeresVel} were
able to establish a zero--full result for simultaneous approximation
with different approximation functions, Theorem BV~\ref{BVsDApcTh1}
below, and a convergence result for multiplicative simultaneous
approximation, Theorem BV~\ref{BVsDApcTh2}, for the class of
non-degenerate rational quadrics.
\begin{theoremBV}\label{BVsDApcTh1}
Let $\psi_1,\psi_2$ be approximating functions and $\mathcal{Q}$ a
$C^{(3)}$ non-degenerate rational quadric. Then
\[
    |\mathcal{Q}\cap\mathcal{S}_2(\psi_1,\psi_2)|_\mathcal{Q}
    = \begin{cases}
    0 & \text{\quad{}if\quad}
    \sum_{q=1}^\infty{}\psi_1(q)\psi_2(q)\leq\infty, \\
    \text{FULL} &\text{\quad{}if\quad}
    \sum_{q=1}^\infty{}\psi_1(q)\psi_2(q) = \infty.
    \end{cases}
\]
\end{theoremBV}
\begin{theoremBV}\label{BVsDApcTh2}
Let $\psi$ be an approximating function and $\mathcal{Q}$ a
$C^{(3)}$ non-degenerate rational quadric. Then
\[
    |\mathcal{Q}\cap\mathcal{S}_2^*(\psi)|_\mathcal{Q} = 0
    \text{\quad{}if\quad}
    \sum_{q=1}^\infty{}\psi_1(q)\psi_2(q) < \infty.
\]
\end{theoremBV}
Of particular note is Theorem BV~\ref{BVsDApcTh2}, which was the
first result of its type for multiplicative simultaneous
approximation on a reasonably general class of manifolds.

The counterpart convergence statement to~\eqref{BDVV1} was finally
established in~\cite{vaugh_velan}. To do so, Vaughan \& Velani
managed to extend Vaughan's result to any sufficiently smooth
function. This result is crucial to our arguments below and we shall
postpone giving the full statement until~\S\ref{ProofThm2} when we
can put the result into a more clearer context.

It should be noted that whilst Beresnevitch \& Velani were able only
to establish the convergence part of Khintchine's theorem for
rational quadrics, they did prove the divergence part in full
generality.

\begin{theoremBV}\label{BVsDApcTh3}
Let $\psi_1,\psi_2$ be approximating functions and $\mathcal{C}$ a
$C^{(3)}$ non-degenerate planar curve. Then
\[
    |\mathcal{C}\cap\mathcal{S}_2(\psi_1,\psi_2)|_\mathcal{C}
    \text{\quad{}is\quad{}FULL\quad{}if\quad}
    \sum_{q=1}^\infty{}\psi_1(q)\psi_2(q) = \infty.
\]
\end{theoremBV}
In \S\ref{ProofThm2} we establish the convergence case for Theorem
BV\ref{BVsDApcTh3} thus completing the Lebesgue metric theory for
this particular problem.

The results discussed above are among the few that are known to hold
for any reasonably general class of manifolds. Obviously there is
still a great deal to be done before the simultaneous theory is
anywhere near as complete as the dual case. A particularly
significant result is due to D.~Kleinbock \& G.~Margulis ~\cite{KM},
who established the validity of the Baker-Sprind\v{z}uk conjecture:
\begin{theoremKM}
\label{KM}
    Let $\mathcal{M}$ be a non-degenerate manifold in $\RR^n$ and
    $\psi:\RR^+\to\RR^+:x\mapsto{}x^{-\tau}$.
    If $\tau>1$ then
    \[
        |\mathcal{M}\cap{}\mathcal{S}^*_n(\psi)|_\mathcal{M} = 0.
    \]
\end{theoremKM}
Theorem KM gives us hope that it will be possible to establish the
analogues of the classical results of Khintchine and Gallagher for
higher dimensional manifolds in much the same way that these results
exist for the dual case. However, there is still a long way to go
before this becomes reality.

\section{Statement of Results}

Let $I$ be some open interval in $\RR$ and $f\in{}C^{(3)}(I)$ such
that for all $x\in{}I$:
\begin{enumerate}
    \item{} there exist constants $c_1>c_2>0$ with $c_1>f'(x)>c_2$,
    \item{} $f^{\prime\prime}\neq0$.
\end{enumerate}
Under these assumptions it is readily verified that the curve $C_f$,
where
\[
    C_f := \{ (x,f(x)):x\in{}I \},
\]
is non-degenerate. We shall be assuming these conditions throughout
the remainder of this article.

In~\cite{vaugh_velan}, the authors conjectured that
\begin{equation}
\label{vv_claim1}
    |\mathcal{C}_f\cap\mathcal{S}_2(\psi_1,\psi_2)|_{\mathcal{C}_f} = 0
    \text{\quad{}if{}\quad} \sum_{h=1}^\infty \psi_1(h)\psi_2(h) <
    \infty
\end{equation}
and
\begin{equation}
\label{vv_claim2}
    |\mathcal{C}_f\cap\mathcal{S}_2^*(\psi)|_{\mathcal{C}_f} = 0
    \text{\quad{}if{}\quad} \sum_{h=1}^\infty \psi(h)\log{}h <
    \infty.
\end{equation}

The conjecture stated in~\eqref{vv_claim2} is a special case of
Theorem~\ref{th2} below. However, before stating the theorem it is
necessary to briefly discuss Hausdorff $s$-measures.

Lebesgue measure is in some sense a relatively ``coarse'' measure of
the size of a set. The notion of measure $0$ can hide a multitude of
finer structure; structure that is often of great interest in
Diophantine approximation. To overcome this technical deficiency,
one can use the idea of Hausdorff $s$-measure and the related notion
of Hausdorff dimension. Hausdorff $s$-measures and dimensions are
(theoretically) computable quantities that give more precise
information about sets of Lebesgue measure $0$ and it is for this
reason that they play such a central role in Diophantine
approximation. We outline only the very basics of the theory of
Hausdorff measures. For a more detailed exposition of the theory and
its many applications in mathematics, see either of the excellent
books by Falconer (~\cite{F1} or~\cite{F2}).

Let $X\subset\RR^n$ and $s\geq0$. For any $\delta>0$, a
$\delta$--cover, $\mathcal{C}_\delta(X)$, of $X$ is a countable
collection of balls $B_i$ such that $ X \subset \bigcup B_i$ and
$\diam{}B_i\leq\delta\}$. The set function
$\mathcal{H}^s_\delta(\cdot)$, where
\[
    \mathcal{H}^s_\delta(X) := \inf \left\{ \sum{}\diam^s{}B_i
    \right\}
\]
with the infimum taken over all $\delta$-covers of $X$, is an outer
measure. Taking the limit of this quantity as $\delta\to{}0$ gives
the Hausdorff $s$-measure of $X$. That is,
\[
    \mathcal{H}^s(X) := \lim_{\delta\to{}0}\mathcal{H}^s_\delta(X).
\]
When $s$ takes on values in $\mathbb{N}\cup\{0\}$ then
$\mathcal{H}^s$ coincides with $s$-dimensional Lebesgue measure. One
of the most useful properties of Hausdorff $s$-measure is the
existence of a unique value of $s$, which we shall denote by
$\dim{}X$, where the $s$-measure jumps from $0$ to $\infty$ as the
 parameter $s$ passes through the value $\dim{}X$ from right to
 left. More precisely,
 \[
    \dim X := \inf \{ s\in\RR^+: \mathcal{H}^s(X) = 0 \} = \sup \{ s\in\RR^+: \mathcal{H}^s(X) = \infty
    \}.
 \]
It is exactly this quantity that we refer to as the Hausdorff
dimension of $X$. At the critical exponent $s=\dim(X)$, the value of
$\mathcal{H}^{\dim(X)}(X)$ can be $0$, $>0$, or $\infty$. Sets $X$
where $\mathcal{H}^s(X) \in (0,\infty)$ are known as
$s$-sets with probably the most famous example of am $s$-set whose Lebesgue measure is $0$ being
 the classical middle-thirds Cantor set
$\mathbb{K}$. Indeed, t is well known that $\dim \mathbb{K} = \ln 2/\ln 3$
and that $\mathcal{H}^{\ln 2/ \ln 3}(\mathbb{K})=1$.  Thus the $s$-measure of a set can be used to garner more precise information
about Lebesgue measure $0$ and this is one of their principle uses in metric Diophantine approximation.

We now come to Theorem~\ref{th2}, which is the Hausdorff $s$-measure
version of Theorem BV~\ref{BVsDApcTh3}.
\begin{theorem}\label{th2}
Let $\psi$ be an approximating function and $0<s \leq 1$. Then
\[
\mathcal{H}^s(C_f\cap S_2^*(\psi))=0\quad\mbox{if}\quad
\sum_{h=1}^{\infty}h^{1-s}(\log^sh)\psi^s(h)<\infty.
\]
\end{theorem}

Note that as $\mathcal{H}^1$ coincides with $1$-dimensional Lebesgue
measure and as we are working with the induced measure on the
manifold, in this case a $1$-dimensional manifold,
Conjecture~\ref{vv_claim2} follows immediately as a special case of
Theorem~\ref{th2}. For the cases when $0 < s < 1$ Theorem~\ref{th2}
appeared as a conjecture in~\cite{vaugh_velan}.

Furthermore, the proof of Theorem~\ref{th2} can be adapted to settle
claim~\ref{vv_claim1} and it is exactly this result that we present
as Theorem~\ref{th1} below.
\begin{theorem}\label{th1}
Let $\psi_1,\psi_2$ be approximating functions. Then
\[
|C_f\cap S_2(\psi_1,\psi_2)|_{C_f} = \begin{cases}
                                            0 \text{\quad{}if\quad}
                                            \sum_{h=1}^{\infty}\psi_1(h)\psi_2(h)<\infty.
                                        \end{cases}
\]
\end{theorem}
Note that as mentioned above, in \S\ref{intro}, this establishes the
convergence counterpart to Theorem BV~\ref{BVsDApcTh3} and this
completes the metric theory, at least in the case of Lebesgue.

\section{Proof of Theorem \ref{th2}}
\label{ProofThm2}

We are given that
\begin{equation}\label{eq25}
\sum_{h=1}^{\infty}h^{1-s}\log^sh\cdot\psi^s(h)<\infty.
\end{equation}
Therefore without loss of generality we can assume that
\begin{equation}\label{eq26}
q^{1-\frac{2}{s}}(\log q)^{-2-\frac{1}{s}}<\psi(q)
\end{equation}
for sufficiently large $q$. To see why, suppose that \eqref{eq26} is
not satisfied. Then we replace $\psi$ with the auxiliary function
\[
\widetilde{\psi}:q \mapsto \widetilde{\psi}:= \max\{\psi(q),
q^{1-\frac{2}{s}}(\log q)^{-2-\frac{1}{s}}\}.
\]
Clearly, $\widetilde{\psi}$ is an approximation function. One can
easily check that~\eqref{eq25} and\eqref{eq26} are satisfied with
$\psi$ replaced by $\widetilde{\psi}$. Furthermore,
\[
S_2^*(\widetilde{\psi})\supset S_2^*(\psi).
\]
Thus it suffices to prove the theorem with $\psi$ replaced by
$\widetilde{\psi}$ and \eqref{eq26} can be assumed.

The set $C_f\cap S_2^*(\psi)$ is a $\limsup$-set with the following
natural representation:
\[
C_f\cap
S_2^*(\psi)=\bigcap_{n=1}^{\infty}\bigcup_{q=n}^{\infty}\bigcup_{(p_1,p_2)\in
\ZZ^2} S^*(p_1,p_2,q)
\]
where
\[
S^*(p_1,p_2,q) := \left\{ (x,y)\in C_f\;:\; \left|
x-\frac{p_1}{q}\right|\cdot\left|
y-\frac{p_2}{q}\right|<\frac{\psi(q)}{q^2}\right\}.
\]

Using the fact that $\psi$ is decreasing, we have that for any $n$
\begin{equation}\label{eq27}
C_f\cap S_2^*(\psi)\subset \bigcup_{t=n}^{\infty}\bigcup_{2^t\le
q<2^{t+1}}\bigcup_{(p_1,p_2)\in \ZZ^2}{}S^*(p_1,p_2,q,t)
\end{equation}
where
\[
    S^*(p_1,p_2,q,t) := \left\{(x,y)\in C_f\;:\;
\left|x-\frac{p_1}{q}\right|\cdot\left|
y-\frac{p_2}{q}\right|<\frac{\psi(2^t)}{(2^t)^2}\right\}.
\]

If  $t\in \NN, (x,y)\in C_f, q\in\NN$ with $2^t\le q<2^{t+1}$ and
\[
\left|x-\frac{p_1}{q}\right|\cdot\left|
y-\frac{p_2}{q}\right|<\frac{\psi(2^t)}{(2^t)^2} \] for some
$(p_1,p_2)\in\ZZ^2$, then there is a unique integer $m$ such that
\[
2^{m-1}\frac{\sqrt{2\psi(2^t)}}{2^t}\le\left|x-\frac{p_1}{q}\right|<2^m\frac{\sqrt{2\psi(2^t)}}{2^t}.
\]
For this number $m$, it follows that
\[
\left|y-\frac{p_2}{q}\right|<2^{-m}\frac{\sqrt{2\psi(2^t)}}{2^t}.
\]
Let
\[
    \gamma_t := \frac{\sqrt{2\psi(2^t)}}{2^t}.
\]
Then, 
\begin{equation}\label{eq28}
C_f\cap S_2^*(\psi)\subset \bigcup_{t=n}^{\infty}\bigcup_{2^t\le
q<2^{t+1}}\bigcup_{(p_1,p_2)\in \ZZ^2}\bigcup_{m=-\infty}^{+\infty}
C_f\cap S(q,p_1,p_2,m)
\end{equation}
where
\[
S(q,p_1,p_2,m)=\left\{(x,y)\in\RR^2\;:\;
\left|x-\frac{p_1}{q}\right|<2^m \gamma_t,
\; \left|y-\frac{p_2}{q}\right|<2^{-m} \gamma_t
\right\}.
\]
Thus, we have constructed a sequence of coverings of $C_f\cap
S_2^*(\psi)$. The aim is now to show that if, for a given $s$,
\eqref{eq25} holds than the associated sequence of Hausdorff
$s$-measures for these coverings tends to 0 as $n\to \infty$. It
then follows that $\mathcal{H}^s(C_f\cap S_2^*(\psi)=0$, as
required.

To proceed we consider two separate cases. For a fixed $t$, Case
(a): $m\in\ZZ$ such that
\begin{equation}\label{eq29}
2^{-|m|}\ge t\sqrt{\psi(2^t)}.
\end{equation}
and Case (b): $m\in\ZZ$ such that
\begin{equation}\label{eq30}
2^{-|m|}\le t\sqrt{\psi(2^t)}.
\end{equation}

\subsection*{Case (a):} First, observe that \eqref{eq29} together
with \eqref{eq26} implies that
\[
2^{-|m|}\ge t\sqrt{2^{t(1-\frac{2}{s})}\cdot
t^{-2-\frac{1}{s}}}\Rightarrow 2^{|m|}\le
t^{\frac{1}{s}}2^{t(\frac{1}{s}-\frac{1}{2})}.
\]
Upon taking logarithms of both sides of the above inequality, we
arrive at
\begin{equation}\label{eq31}
|m|\le t\left(\frac{2-s}{2s}\right) +\frac{1}{2s}\log t\ll t.
\end{equation}
As  $f'(x)>c_1$ for all $x\in I$ it follows that
\begin{equation}\label{eq32}
\diam(C_f\cap S(q,p_1,p_2,m)) \ll
2^{-|m|}\frac{\sqrt{\psi(2^t)}}{2^t}.
\end{equation}
The implied constant depends on only $c_1$ and is 
irrelevant to the remainder of the argument.

Given $t$ and $m$, let $N(t,m)$ denote the number of triples
$(q,p_1,p_2)$ with $2^t\le q<2^{t+1}$ such that $C_f\cap
S(q,p_1,p_2,m)\neq \emptyset$. Suppose now that $C_f\cap
S(q,p_1,p_2,m)\neq \emptyset$. Then for some $(x,y)\in C_f$ and
$\theta_1,\theta_2$ satisfying $-1<\theta_1,\theta_2<1$, we have
that
\[
x=\frac{p_1}{q}+\theta_1 2^{|m|}\frac{\sqrt{2\psi(2^t)}}{2^t},\quad
y=\frac{p_2}{q}+\theta_2 2^{|m|}\frac{\sqrt{2\psi(2^t)}}{2^t}.
\]
Thus, it can be shown that
\begin{align*}
f\left(\frac{p_1}{q}\right)-\frac{p_2}{q}
&=f\left(\frac{p_1}{q}\right)-f(x)+f(x)-y+y-\frac{p_2}{q}\\
&= -\theta_2f'(\xi)\cdot
2^{|m|}\frac{\sqrt{2\psi(2^t)}}{2^t}+\theta_1\cdot
2^{|m|}\frac{\sqrt{2\psi(2^t)}}{2^t}
\end{align*}
where $\xi$ lies between $x$ and $p_1/q$. Further, one can easily
deduce that
\[
\left|f\left(\frac{p_1}{q}\right)-\frac{p_2}{q}\right|\ll
2^{|m|}\frac{\sqrt{\psi(2^t)}}{2^t}\le \frac{1}{t2^t}.
\]

Set $Q=2^{t+1}$. Then we have $q\le Q$ and
\[
\left|f\left(\frac{p_1}{q}\right)-\frac{p_2}{q}\right|\ll
\frac{1}{Q\log Q}.
\]
As mentioned in the introductory section of this article, a result
of Vaughan \& Velani is crucial to our argument. We now state this
result, which is Theorem 2 from \cite{vaugh_velan}.
\begin{theoremVV}
Let $N_f(Q,\psi,I)=\#\{\mathbf{p}/q\::\: q\le Q, p_1/q\in I,
|f(p_1/q)-p_2/q|<\psi(Q)/Q\}$. Suppose that $\psi$ is an
approximating function with $\psi(Q)\ge Q^{-\phi}$ where $\phi$ is
any real number with $\phi\le \frac{2}{3}$. Then
\begin{equation}
N_f(q,\psi,I)\ll \psi(Q)Q^2.
\end{equation}
\end{theoremVV}
In our case $\psi(Q)=2^{|m|}\sqrt{2\psi(2^t)}\asymp \frac{1}{\log
Q}$ which satisfies the conditions of Theorem VV. Therefore there
exists an absolute constant $c>0$ such that
\begin{equation}\label{eq34}
N(t,m)\ll 2^{2t}2^{|m|}\sqrt{\psi(2^t)}.
\end{equation}

Now using~\eqref{eq32} and~\eqref{eq34} we can bound the Hausdorff
sum 
associated with the set $C_f\cap S_2^*(\psi)$
\begin{align*}
\mathcal{H}^s(C_f\cap S_2^*(\psi)) &\ll
\sum_{t=n}^{\infty}\sum_{m\in\mbox{case
(a)}}\left(2^{-|m|}\frac{\sqrt{\psi(2^t)}}{2^t}\right)^s\times
2^{|m|}2^{2t}\sqrt{\psi(2^t)} \\
&\ll \sum_{t=n}^\infty\sum_{m\in\mbox{case (a)}}
\psi(2^t)^{\frac{1}{2}(1+s)}\cdot 2^{t(2-s)}\cdot 2^{|m|(1-s)} \\
&\stackrel{(\ref{eq28})}{\ll}\sum_{t=n}^\infty\sum_{m\in\mbox{case
(a)}} t^{s-1}\psi(2^t)^s2^{t(2-s)}\\
&\stackrel{(\ref{eq31})}{\ll} \sum_{t=n}^{\infty}
t^s2^{t(2-s)}\psi(2^t)^s\asymp \sum_{q=2^n}^{\infty} q^{1-s}\log^s
q\cdot\psi(q)^s.
\end{align*}
The above comparability follows from the fact that $\psi$ is an
approximating function and therefore decreasing. In view
of~\eqref{eq25},
\[
\sum_{q=2^n}^{\infty} q^{1-s}\log^s q\cdot\psi^s(q)\to 0\;\mbox{ as
}\; n\to \infty.
\]
Therefore, for Case (a) it follows that $\mathcal{H}^s(C_f\cap
S_2^*(\psi))=0$ as required.

\subsection*{Case (b):}

In view of~\eqref{eq30}, we have that
\[
S(q,p_1,p_2,m)\subset \left(S^\prime(q,p_1)\times [0,1]\right)\cup
\left([0,1]\times S^\prime(q,p_2)\right)
\]
where
\[
S'(q,p)=\left\{y\in [0,1]\::\:
\left|y-\frac{p}{q}\right|<\frac{2t\psi(2^t)}{2^t}\right\}.
\]

Thus, the set of~\eqref{eq28} 
is a subset of
\begin{equation}\label{eq37}
\bigcup_{t=n}^{\infty}\bigcup_{2^t\le
q<2^{t+1}}\bigcup_{(p_1,p_2)\in \ZZ^2} C_f\cap
\left(\left(S^\prime(q,p_1)\times[0,1]\right)\cup \left([0,1]\times
S^\prime(q,p_2)\right)\right).
\end{equation}

As $c_1>f'(x)>c_2>0$, for any choice of $p_1,p_2$ and $q$ appearing
in~\eqref{eq37},
\[
\diam(C_f\cap(S'(q,p_1)\times [0,1]))\ll
t\frac{\psi(2^t)}{2^t}
\]
and
\[
\diam(C_f\cap([0,1]\times S'(q,p_2))\ll t\frac{\psi(2^t)}{2^t}.
\]
The implied constants depends only on $c_1$ and $c_2$ and are
irrelevant in the context of the rest of the proof. Furthermore, for
a fixed $t$ and $q$ in~\eqref{eq37}, the number of
$(p_1,p_2)\in\ZZ^2$ for which
\[
C_f\cap \left(\left(S'(q,p_1)\times [0,1]\right)\cup
\left([0,1]\times
S'(q,p_2)\right)\right)
\]
are non-empty and disjoint is $\ll q$. Drawing all the above
considerations together it follows that the Hausdorff $s$-sum for
this covering of the set $C_f\cap S_2^*(\psi)$, as defined
in~\eqref{eq37}, is bounded above by
\[
\sum_{t=n}^{\infty}\left(t\frac{\psi(2^t)}{2^t}\right)^s2^{2t}
\asymp \sum_{q=2^n}^{\infty} q^{1-s}\log^sq\cdot\psi^s(q)\to
0\;\mbox{ as } n\to \infty.
\]
As in the previous case, Case (a), by letting $n\to\infty$ we
conclude that
\[
    \mathcal{H}^s(C_f\cap S_2^*(\psi))=0
\]

This completes the proof of Theorem~\ref{th2}.

\section{Proof of Theorem~\ref{th1}}

We now proceed with establishing Theorem~\ref{th1}. The proof is
somewhat analogous to that of Theorem~\ref{th2} and for brevity we
leave the technical details needed to verify some of our estimates
to the reader. Especially in those cases when the estimates
following in exactly the same manner or require only minor
modifications to the arguments used in the proof Theorem~\ref{th2}.

For the sake of convenience, let $\psi=\psi_1$ and $\phi=\psi_2$. It
is clear that
\[
S_2(\phi,\psi)\subset S_2(\psi^*,\psi_*)\cup S_2(\psi_*,\psi^*)
\]
where
\[
\psi_*=\min\{\psi,\phi\} \text{\quad\&\quad}
\psi^*=\max\{\psi,\phi\}.
\]
Since $\psi^*\psi_*=\psi\phi$, we have that
$\sum\psi^*(q)\psi_*(q)<\infty$. Thus to prove Theorem~\ref{th1} it
is sufficient to show that both the sets $C_f\cap
S_2(\psi^*,\psi_*)$ and $C_f\cap S_2(\psi_*,\psi^*)$ are of Lebesgue
measure zero. We will consider one of these two sets, the other case
is similar. Thus, without any loss of generality we assume that $
\psi(q)\ge \phi(q)
$
for all $q\in\NN$.

Since $\sum_{q=1}^\infty\psi(q)\phi(q)<\infty$ and both $\psi,\phi$
are decreasing we have that $\psi(q)\phi(q)<q^{-1}$ for all
sufficiently large $q$. Hence, $ \phi(q)\le q^{-1/2}$ for
sufficiently large $q$. Further, we can assume that
\begin{equation}\label{eq38}
\psi(q)\ge q^{-2/3}
\end{equation}
for all $q\in \mathbb{N}$. To see this consider the auxiliary
function $\tilde{\psi}$ where
 \[
 \tilde{\psi}(q)=\max\{\psi(q),q^{-2/3}\}.
 \]
Clearly, $\tilde{\psi}$ is an approximating function. It also
satisfies the following set inclusion,
\[
S_2(\psi,\phi)\subset S_2(\tilde{\psi},\phi).
\]
Moreover,
\begin{align*}
\sum_{q=1}^\infty \tilde{\psi}(q)\phi(q) & \le
\sum_{q=1}^\infty\psi(q)\phi(q)+\sum_{q=1}^\infty q^{-2/3}\phi(q)
\\
&\ll \sum_{q=1}^\infty \psi(q)\phi(q)+\sum_{q=1}^\infty
q^{-2/3}q^{-1/2}<\infty.
\end{align*}
This means that it is sufficient to prove Theorem~\ref{th1} with
$\psi$ replaced by $\tilde{\psi}$ and therefore without any loss of
generality,~\eqref{eq38} can be assumed.

In a manner analogous to that of~\eqref{eq27}, it is readily
verified that for any $n\ge 1$
\begin{equation}
\label{th2cover}
C_f\cap S_2(\psi,\phi)\subset
\bigcup_{t=n}^\infty\bigcup_{2^t\le q<2^{t+1}}\bigcup_{(p_1,p_2)\in
\mathbb{Z}^2}C_f\cap S_2(p_1,p_2,q)
\end{equation}
where
\[
S_2(p_1,p_2,q)=\left\{ (x,y)\in \mathbb{R}^2 : \left|
x-\frac{p_1}{q}\right|<\frac{\psi(2^t)}{2^t}, \left|
y-\frac{p_2}{q}\right|<\frac{\phi(2^t)}{2^t}\right\}
\]
and $t$ is uniquely defined by $2^t\le q<2^{t+1}$. Next, we can use
the same argument to that used in~\eqref{eq32} to verify that
\begin{equation}
\label{th2coverest}
|C_f\cap S_2(q,p_1,p_2)|_{C_f}\ll
\frac{\phi(2^t)}{2^t}.
\end{equation}
Finally, for fixed $t$ let $N(t)$ be the number of triples
$(q,p_1,p_2)$ with $2^t\le q<2^{t+1}$ such that $C_f\cap
S(q,p_1,p_2)\neq \emptyset$. On modifying the argument used to
establish~\eqref{eq34}, one obtains the estimate
\begin{equation}
\label{th2Ntest}
N(t)\ll 2^{2t}\psi(2^t).
\end{equation}

The upshot of~\eqref{th2cover},~\eqref{th2coverest} and~\eqref{th2Ntest} is that
\begin{align*}
|C_f\cup S_2(\psi,\phi)|_{C_f} & \ll \sum_{t=n}^\infty \sum_{2^t\le
q<2^{t+1}} \bigcup_{(p_1,p_2)\in \mathbb{Z}^2}C_f\cap
S_2(q,p_1,p_2)\\
&\ll \sum_{t=n}^\infty N(t)\frac{\phi(2^t)}{2^t}\ll
\sum_{t=n}^\infty 2^t\psi(2^t)\phi(2^t) \asymp \sum_{q=2^n}^\infty
\psi(q)\phi(q).
\end{align*}
Since 
$\sum_{q=1}^\infty \psi(q)\phi(q)<\infty$, we have that
$\sum_{q=2^n}^\infty \psi(q)\phi(q)\to 0$ as $n\to\infty$ and it
follows that
\[
|C_f\cap S_2(\psi,\phi)|_{C_f}=0
\]
as required.

This completes the proof of Theorem~\ref{th1}.

\section{Remarks and Possible Developments}

An obvious next step is to establish the divergence counterpart to
Theorem~\ref{th2}. That is, one would like to show that
\[
    \sum_{h=1}^\infty h^{1-s}(\log^sh)\psi^s(h) = \infty
    \implies
    \mathcal{H}^s(C_f\cap S_2^*(\psi))= \infty.
\]
By adapting the arguments in this paper and using the ideas of local
ubiquity, as developed in~\cite{BDV0401118}, it is likely that one
could establish a zero-full result for $\mathcal{H}^h(C_f\cap
S_2^*(\psi))$ where $h$ is a general dimension function. This would
include the above result and Theorem~\ref{th2} as a special case.  A
dimension function $h:\RR^+\to\RR^+$ is an increasing, continuous
function such that $h(r)\to0$ as $r\to0$. By replacing the quantity
$\diam^s(C_i)$ with $h(\diam(C_i))$ in the definition of
$\mathcal{H}^s$, one can define the Hausdorff $h$-measure of a set.
For further details see~\cite{F2} or~\cite{Mat}. Dimension functions
give very precise information about the measure theoretic properties
of a set. The convergence part of such a theorem follows almost
immediately on from Theorem~\ref{th2}. Most of the estimates
obtained in the proof of Theorem~\ref{th2} remain the same, the
generalisation to Hausdorff $h$-measures effects only the estimates
involving the measures of the actual covers defined in the proof.
The main task in proving such a theorem would be in the proof of the
divergence case.

\textbf{Acknowledgements}~JL thanks ``Mr. Diplomacy'' himself, Sanju
Velani, for his constant help and encouragement in all matters
mathematical (and sometimes political) and for keeping me on the
straight and narrow.

\end{document}